\documentclass[12pt,reqno]{amsart}
\usepackage{amsmath}
\usepackage{amsxtra}
\usepackage{amscd}
\usepackage{amsthm}
\usepackage{amsfonts}
\usepackage{amssymb}
\usepackage{eucal}

\newtheorem{theorem}{Theorem}[section]
\newtheorem{conj}[theorem]{Conjecture}
\newtheorem{cor}[theorem]{Corollary}
\newtheorem{lem}[theorem]{Lemma}

\theoremstyle{definition}

\theoremstyle{remark}
\newtheorem{rem}[theorem]{Remark}
\theoremstyle{remark}

\numberwithin{equation}{section}

\def\slg{\mathfrak{sl}}
\def\glg{\mathfrak{gl}}

\let\leq\leqslant

\let\geq\geqslant

\newcommand{\nc}{\newcommand}
\nc{\on}{\operatorname}
\nc{\ch}{\mbox{ch}}
\nc{\Z}{{\mathbb Z}}
\nc{\C}{{\mathbb C}}
\nc{\R}{{\mathbb R}}
\nc{\pone}{{\mathbb C}{\mathbb P}^1}
\nc{\pa}{\partial}
\nc{\F}{{\mathcal F}}
\nc{\arr}{\rightarrow}
\nc{\larr}{\longrightarrow}
\nc{\al}{\alpha}
\nc{\ri}{\rangle}
\nc{\lef}{\langle}
\nc{\W}{{\mathcal W}}
\nc{\la}{\lambda}
\nc{\ep}{\epsilon}
\nc{\eps}{\varepsilon}

\nc{\Om}{\Omega}
\nc{\su}{\widehat{{\mathfrak sl}}_2}
\nc{\sw}{{\mathfrak s}{\mathfrak l}}

\nc{\g}{{\mathfrak g}}
\nc{\h}{{\mathfrak h}}
\nc{\n}{{\mathfrak n}}
\nc{\N}{\widehat{\n}}
\nc{\G}{\widehat{\g}}
\nc{\De}{\Delta_+}
\nc{\gt}{\widetilde{\g}}
\nc{\Ga}{\Gamma}
\nc{\one}{{\mathbf 1}}
\nc{\z}{{\mathfrak Z}}
\nc{\zz}{{\mathcal Z}}
\nc{\Hh}{{\mathcal H}_\beta}
\nc{\qp}{q^{\frac{k}{2}}}
\nc{\qm}{q^{-\frac{k}{2}}}
\nc{\La}{\Lambda}
\nc{\wt}{\widetilde}
\nc{\qn}{\frac{[m]_q^2}{[2m]_q}}
\nc{\cri}{_{\on{cr}}}
\nc{\kk}{h^\vee}
\nc{\sun}{\widehat{\sw}_N}
\nc{\hh}{\widehat{\mathfrak h}}
\nc{\HH}{{\mathcal H}_{q,t}}
\nc{\ca}{\wt{{\mathcal A}}_{h,k}(\sw_2)}
\nc{\gl}{\widehat{{\mathfrak g}{\mathfrak l}}_2}
\nc{\el}{\ell}
\nc{\s}{{\mathbf s}}
\nc{\bi}{\bibitem}
\nc{\om}{\omega}
\nc{\WW}{\W_\beta}
\nc{\scr}{{\mathbf S}}
\nc{\ab}{{\mathbf a}}
\nc{\rr}{r}
\nc{\ol}{\overline}
\nc{\con}{qt^{-1} + q^{-1}t}
\nc{\den}{q^{\el-1} t^{-\el+1}+ q^{-\el+1} t^{\el-1}}
\nc{\ds}{\displaystyle}
\nc{\B}{B}
\nc{\A}{{\mathbb A}}
\nc{\GG}{{\mathcal G}}
\nc{\UU}{{\mathcal U}}
\nc{\MM}{{\mathcal M}}
\nc{\CC}{{\mathcal C}}
\nc{\GL}{{}^L G}
\nc{\dzz}{\frac{dz}{z}}
\nc{\Res}{\on{Res}}
\nc{\rep}{{\mathcal R}ep \;}
\nc{\uqg}{U_q \G}
\nc{\uqgg}{U_q \g}
\nc{\Fq}{{\mathbb F}_q}

\nc{\stimes}{\ltimes}
\nc{\K}{\hat{\mathcal K}}
\nc{\Ql}{\ol{\mathbb Q}_\ell}

\nc{\ga}{\gamma}
\nc{\PL}{{}^L P}
\nc{\E}{\mc E}
\nc{\mc}{\mathcal}
\nc{\mbf}{\mathbf}
\nc{\bb}{{\mathfrak b}}
\nc{\OO}{{\mc O}}
\nc{\Po}{{\mc P}}
\nc{\V}{{\mc V}}
\nc{\yy}{{\mc Y}}
\nc{\M}{\mathcal M}
\nc{\Coh}{{{\mathcal C}oh}}
\nc{\Cohn}{\Coh_n}
\nc{\f}{{\mathcal F}}
\nc{\si}{_E}
\nc{\Gaf}{{\mathbb G}_{a,\Fq}}
\nc{\KK}{{\mathfrak k}}

\nc{\PO}{{\mathbb P^1}}
\nc{\PR}{{\mathbb P^r}}

\nc{\Wr}{{ {\rm Wr}}}

\newcommand{\bean}{\begin{eqnarray}}
\newcommand{\eean}{\end{eqnarray}}
\newcommand{\be}{\begin{displaymath}}
\newcommand{\ee}{\end{displaymath}}
\newcommand{\bea}{\begin{eqnarray*}}
\newcommand{\eea}{\end{eqnarray*}}
\newcommand{\beq}{\begin{equation}}
\newcommand{\eeq}{\end{equation}}
\newcommand{\bs}{\boldsymbol}
\newcommand{\Ref}[1]{{$($\ref{#1}$)$}}

\newcommand{\sing}{{\rm Sing}\,}

\newcommand{\nL}{L_{\om_r}^{\otimes n}[\mu]}
\newcommand{\nnL}{L_{\om_r}^{\otimes n}}

\newcommand{\snL}{ \sing L_{\om_r}^{\otimes n}[\mu]}
\newcommand{\btz}{ \om(\bs t; \bs z)}

\newcommand{\GR}{ {G(r+1,d)}}

\begin{document}

\title[The B. and M. Shapiro conjecture and the Bethe ansatz]
{The B. and M. Shapiro conjecture in real algebraic geometry
and the Bethe ansatz}

\author[E. Mukhin, V. Tarasov, and A. Varchenko]
{E. Mukhin ${}^{*}$, V. Tarasov ${}^{*,\star,1}$,
\and A. Varchenko {${}^{**,2}$} }
\thanks{${}^1$\ Supported in part by RFFI grant 05-01-00922}
\thanks{${}^2$\ Supported in part by NSF grant DMS-0244579}

\maketitle

\centerline{\it ${}^*$Department of Mathematical Sciences,
Indiana University -- Purdue University,}
\centerline{\it Indianapolis, 402 North Blackford St, Indianapolis,
IN 46202-3216, USA}
\smallskip
\centerline{\it $^\star$St.\,Petersburg Branch of Steklov Mathematical
Institute}
\centerline{\it Fontanka 27, St.\,Petersburg, 191023, Russia}
\smallskip
\centerline{\it ${}^{**}$Department of Mathematics, University of
North Carolina at Chapel Hill,} \centerline{\it Chapel Hill, NC
27599-3250, USA} \medskip

\medskip

\begin{abstract}
We prove the B. and M. Shapiro conjecture: if the Wronskian of
a set of polynomials has real roots only, then the complex span of this set
of polynomials has a basis consisting of polynomials with real
coefficients. This, in particular, implies the following result:

If all ramification points of a parametrized rational curve
$\phi:\C\mathbb P^1 \to
\C\mathbb P^r$ lie on a circle in the Riemann sphere $\C\mathbb P^1$,
then $\phi$ maps this circle into a suitable real subspace
$\mathbb R\mathbb P^r \subset \C\mathbb P^r$.

The proof is based on the Bethe ansatz method in the Gaudin model. The
key observation is that a symmetric linear operator on a Euclidean
space has real spectrum.

In Appendix A, we discuss properties of differential operators
associated with Bethe vectors in the Gaudin model. In particular, we
prove a statement which may be useful in complex algebraic geometry:
certain Schubert cycles in a Grassmannian intersect transversally if
the spectrum of suitable Gaudin Hamiltonians is simple.

In Appendix B, we formulate a conjecture on
reality of orbits of critical points of master functions and prove this
conjecture for master functions associated with Lie algebras of types
$A_r, B_r, C_r$.
\end{abstract}

\maketitle

\thispagestyle{empty}

\section{The B. and M. Shapiro conjecture}

\subsection{Statement of the result}\label{statement section}
Fix a natural number $r\geq 1$. Let $V\subset\C[x]$ be a vector subspace of
dimension $r+1$. The space $V$ is called {\it real} if it has a basis
consisting of polynomials in $\R[x]$.

For a given $V$, there exists a unique linear differential operator
$$
D\ =\ \frac {d^{r+1}}{dx^{r+1}} \ +\ \la_1(x)\, \frac {d^{r}}{dx^{r}}\ +
\dots +\ \la_{r}(x)\, \frac {d^{}}{dx^{}}\ +\ \la_{r+1}(x)\ ,
$$
whose kernel is $V$. This operator is called {\it the fundamental differential
operator of $V$}. The coefficients of the operator are rational functions in $x$.
The space $V$ is real if and only if all coefficients of the fundamental operator are
real rational functions.

{\it The Wronskian} of functions $f_1,\dots,f_i$ in $x$ is the
determinant
\[
\Wr (f_1,\dots, f_i)\ =\ \det
\left( \begin{array} {cccccc}
f_1 & f_1^{(1)} & \dots\dots & f_1^{(i-1)}
\\
f_2 & f_2^{(1)} &{} \dots\dots & f_2^{(i-1)}
\\
{}\dots & {} \dots & \dots\dots & \dots
\\
f_i & f_i^{(1)} &{} \dots\dots & f_i^{(i-1)}
\end{array} \right).
\]

Let $f_1,\dots, f_{r+1}$ be a basis of $V$. The Wronskian of the basis
does not depend on the choice of the basis up to multiplication by a
number. The monic representative is called {\it the Wronskian of $V$}
and denoted by ${\rm Wr}_V$.

\begin{theorem}
\label{theorem 1}
If all roots of the polynomial ${\rm Wr}_V$ are real, then the space $V$ is real.
\end{theorem}

This statement is the B. and M. Shapiro conjecture formulated in 1993. 
The conjecture is
proved in \cite{EG1} for $r=1$, see a more elementary proof also
for $r=1$ in \cite{EG3}.
The conjecture, its supporting
evidence, and applications
are discussed in \cite{EG1}\,--\,\cite{EG3}, \cite{EGSV}, \cite{ESS},
\cite{KS}, \cite{RSSS}, \cite{S1}\,--\,\cite{S6}.

\subsection{Parametrized rational curves with real ramification points}
For a projective coordinate system $(v_1:\dots :v_{r+1})$ on the
complex projective space $\C\PR$, the subset of points with real
coordinates is called {\it the real projective subspace} and denoted by
$\mathbb R\PR$.

Let $\phi : \C\PO \to \C\PR$ be a parametrized rational curve.
If $(u_1:u_2)$ are projective coordinate on $\C\PO$ and
$(v_1:\dots :v_{r+1})$ are projective coordinate on $\C\PR$,
then $\phi$ is given by the formula
\bea
\phi\ :\ (u_1:u_2) \ \mapsto\
(\phi_1(u_1,u_2) : \dots : \phi_{r+1}(u_1,u_2))
\eea
where $\phi_i$ are homogeneous polynomials of the same degree.
We assume that at any point of $\C\PO$ at least one of $\phi_i$ is nonzero.
Choose the local affine coordinate $u = u_1/u_2$ on $\C\PO$
and local affine coordinates
$v_1/v_{r+1},\dots,v_r/v_{r+1}$ on $\C\PR$. In this coordinates, the map $\phi$ takes the
form
\bean
\label{loc coord}
f\ :\ u \ \mapsto \ \left(\ \frac{f_1(u)}{f_{r+1}(u)}\ ,\ \dots \ ,\
\frac{f_r(u)}{f_{r+1}(u)}\ \right)
\eean
where $f_i(u) = \phi_i(u,1)$.

The map $\phi$ is said
to be {\it ramified} at a point of $\C\PO$ if its first $r$ derivatives at this
point do not span $\C\PR$~\cite{KS}. More precisely, a point $u$ is {\it
a ramification point}, if
the vectors $f^{(1)}(u), \dots, f^{(r)}(u)$ are linear dependent.

We assume that a generic point of $\C\PO$
is not a ramification point.

\begin{theorem}
\label{theorem on ram points}
If all ramification points of the
parametrized rational curve $\phi$ lie on a circle in the Riemann
sphere $\C\mathbb P^1$, then $\phi$ maps this circle into a suitable
real subspace $\mathbb R\mathbb P^r \subset \C\mathbb P^r$.
\end{theorem}

A {\it maximally inflected curve} is, by definition \cite{KS}, a
parametrized real rational curve, all of whose ramification points are real.
From Theorem \ref{theorem on ram points}, it follows the existence of
maximally inflected curves, for every placement of the ramification points.

\medskip

Theorem \ref{theorem on ram points} follows from Theorem \ref{theorem
1}. Indeed, if all ramification points lie on a circle, then changing
linearly the coordinates $(u_1:u_2)$, we may assume that the
ramification points lie on the real line $\R\PO$ and the point $(0:1)$
is not a ramification point. Changing linearly the coordinates
$(v_1:\dots :v_{r+1})$ on $\C\PR$, we may assume that $\phi_{r+1}$ is
not zero at any of the ramification points. Let us use the affine
coordinates $u=u_1/u_2$ and $v_1/v_{r+1}, \dots , v_r/v_{r+1}$, and
formula \Ref{loc coord}. Then the determinant of coordinates of the
vectors $f^{(1)}(u), \dots, f^{(r)}(u)$ is equal to
\bea
\Wr \left(\frac{f_1}{f_{r+1}}, \dots , \frac{f_r}{f_{r+1}}, 1 \right)\!\!(u)
=
\frac{1}{(f_{r+1})^{r+1}}
\Wr (f_1, \dots , f_r, f_{r+1})(u) \ .
\eea
Hence the vectors $f^{(1)}(u), \dots, f^{(r)}(u)$ are linearly
dependent if and only if the Wronskian of $f_1,\dots,f_{r+1}$ at $u$
is zero. Since not all points of $\C\PO$ are ramification points, the
complex span $V$ of polynomials $f_1,\dots,f_{r+1}$ is an
$r+1$-dimensional space. By assumptions of Theorem \ref{theorem on ram
points}, all zeros of the Wronskian of $V$ are real. By Theorem
\ref{theorem 1}, the space $V$ is real. This means that there exist
projective coordinates on $\C\PR$, in which all polynomials
$f_1,\dots,f_{r+1}$ are real. Theorem \ref{theorem on ram points} is
deduced from Theorem \ref{theorem 1}.

\subsection{Reduction of Theorem \ref{theorem 1} to a special case}

\begin{theorem}
\label{theorem 2}
Assume that all roots of the Wronskian are real and simple, then
$V$ is real.
\end{theorem}

We deduce Theorem \ref{theorem 1} from Theorem \ref{theorem 2}.
Indeed, let $V_0$ be an $r+1$-dimensional space of polynomials whose Wronskian
has real roots only. Let $d$ be the degree of a generic polynomial in
$V_0$. Denote
\begin{enumerate}
\item[$\bullet$]
$\C_d[x]$ the space of polynomials of degree not
greater than $d$,
\item[$\bullet$]
$G(r+1,d)$ the Grassmannian of
$r+1$-dimensional vector subspaces in $\C_d[x]$,
\item[$\bullet$]
$\mathbb P(\C_{(r+1)(d-r)}[x])$ the projective space
associated with the vector space
$\C_{(r+1)(d-r)}[x]$.
\end{enumerate}
The varieties $\GR$ and $\mathbb P(\C_{(r+1)(d-r)}[x])$
have the same dimension. The assignment $V \mapsto \Wr_V$ defines a
finite morphism $\pi : \GR \to \mathbb P(\C_{(r+1)(d-r)}[x])$, see,
for example, \cite{S2, EG1}. The space $V_0$ is a point of $\GR$.

Since $\pi$ is finite and $V_0$ has Wronskian with real roots only,
there exists a continuous curve $\ep\mapsto V_\ep\in\GR$ for $\ep\in[0,1)$,
such that the Wronskian of $V_{\ep}$ for $\ep > 0$ has simple real roots only.
By Theorem \ref{theorem 2}, the space $V_\ep$ is real for $\ep > 0$. Hence, the
fundamental differential operator of $V_\ep$ has real coefficients.
Therefore, the fundamental differential operator of $V_0$ has real coefficients
and the space $V_0$ is real. Theorem \ref{theorem 1} is deduced from
Theorem~\ref{theorem 2}.

\subsection{The upper bound for the number of complex
vector spaces with the same exponents
at infinity and the same Wronskian}
\label{number of spaces section}
Let $f_1,\dots,f_{r+1}$ be a basis of $V$ such that $\deg f_i = d_i$ for
some sequence
$$
\bs d\ =\ \{\,d_1\,<\,\dots\,< \,d_{r+1}\,\}\ .
$$
We say that $V$ has exponents $\bs d$
at infinity. If $V$ has exponents $\bs d$ at
infinity, then $\deg \Wr_V = n$ where
$$
n \,=\,
\sum_{i=1}^{r+1}\,(d_i-i+1)\ .
$$
Let
$$
T\ =\ \prod_{s=1}^n\,(x-z_s)
$$
be a polynomial with simple (complex) roots.
Then the upper bound for the number of complex
vector spaces $V$ with exponents $\bs d$ at
infinity and Wronskian $T$
is given by the number $N(\bs d)$ defined as follows.

Consider the Lie algebra $\slg_{r+1}$ with Cartan decomposition
$\slg_{r+1}=\n_-\oplus \h\oplus \n_+$ and simple roots
$\al_1,\dots,\al_{r} \in \h^*$. Fix the invariant inner product
on $\h^*$ by the condition $(\al_i,\al_i)=2$. For any integral dominant weight
$\La \in \h^*$, denote by $L_\La$ the irreducible $\slg_{r+1}$-module with
highest weight $\La$. Let $\om_r \in \h^*$ be the last fundamental weight.

For $i=1,\dots, r$, introduce the numbers
\bean
\label{def of l_i}
l_i \ =\ \sum_{j=1}^i\ ( d_j-j+1)\ ,
\eean
and the integral dominant weight
\bea
\La(\bs d) \ =\ n\om_{r}\ -\ \sum_{i=1}^{r}\,l_i\,\al_i\ .
\eea
Set $N(\bs d)$ to be the multiplicity of the module
$L_{\La(\bs d)}$ in the $n$-factor tensor product
\bea
L_{\om_{r}}^{\otimes n}\ = L_{\om_r} \otimes \dots \otimes L_{\om_r}\ .
\eea

According to Schubert calculus, the number of complex
$r+1$-dimensional vector spaces $V$ with exponents $\bs d$ at
infinity and Wronskian $T$ is not greater than the number
$N(\bs d)$.

This is a standard Schubert calculus statement,
see, for example,
\cite[Section~5]{MV2}.

Thus, in order to prove Theorem \ref{theorem 2}, it is enough to prove

\begin{theorem}
\label{theorem 3}
For generic real $z_1, \dots, z_n$, there exists exactly $N(\bs d)$ distinct
real vector spaces $V$ with exponents $\bs d$ at infinity and Wronskian
$T = \prod_{s=1}^n(x-z_s)$.
\end{theorem}

\subsection{Structure of the paper}
\label{Structure of the paper}
In Section \ref{Construction of spaces of polynomials},
for generic complex $z_1, \dots , z_n$,
we will construct exactly $N(\bs d)$ distinct complex vector spaces
$V$ with exponents $\bs d$ at infinity and Wronskian $T$.
In Section \ref{Bethe vectors}, we will show
that all of these vector spaces are real,
if $z_1,\dots,z_n$ are real. This will prove Theorem \ref{theorem 3}.

The constructions of Sections \ref{Construction of spaces of polynomials}
and \ref{Bethe vectors} are the Bethe ansatz constructions for the Gaudin model
on $L_{\om_{r}}^{\otimes n}$.

In Appendix A, we discuss properties of differential operators associated with
the Bethe vectors in the Gaudin model and give applications of the Bethe ansatz
constructions of Section~\ref{Bethe vectors}.
In particular, we
prove a statement which may be useful in complex algebraic geometry:
certain Schubert cycles in a Grassmannian intersect transversally if
the spectrum of suitable Gaudin Hamiltonians is simple, see
Corollary \ref{last cor}, cf.~\cite{EH} and \cite{MV2}.

In Appendix B, we formulate a conjecture on
reality of orbits of critical points of master functions and prove this
conjecture for master functions associated with Lie algebras of types
$A_r, B_r, C_r$.

\medskip

We thank A.~Eremenko and A.~Gabrielov for useful discussions.

\section{Construction of spaces of polynomials}
\label{Construction of spaces of polynomials}

\subsection{Construction of (not necessarily real)
spaces with exponents $\bs d$ at infinity and
Wronskian $T = \prod_{s=1}^n(x-z_s)$ with simple roots}
\label{construction of spaces section}
Denote $\bs z = (z_1,\dots,z_n)$.
Introduce a function of $l_1+\dots+l_{r}$ variables
\bea
&
\bs t = (t^{(1)}_{1},\dots,t_{l_1}^{(1)},\dots,
t^{(r)}_{1},\dots,t_{l_{r}}^{(r)})\
\eea
by the formula
\bean
\label{master}
\Phi_{\bs d} (\bs t;\bs z) =
\prod_{j=1}^{l_{r}}\prod_{s=1}^n
(t_j^{(r)}-z_s)^{-1}
\prod_{i=1}^{r}\prod_{1\leq j<s\leq l_i} (t_j^{(i)}-t_s^{(i)})^{2}
\prod_{i=1}^{r-1}\prod_{j=1}^{l_i}\prod_{k=1}^{l_{i+1}}
(t_j^{(i)}-t_k^{(i+1)})^{-1} \ .
\eean
The function $\Phi_{\bs d}$ is a rational function of $\bs t$, depending on
parameters $\bs z$. The function is called {\it the master function}.

The master functions arise in the hypergeometric solutions
of the KZ equations
\cite{SV, V1} and in the Bethe ansatz method for the Gaudin
model \cite{RV, ScV, MV1, MV2, MV3, V2}.

The product of symmetric groups $\Sigma_{\bs l}=\Sigma_{l_1}\times \dots \times
\Sigma_{l_r}$ acts on the variables $\bs t$ by permuting the coordinates with the
same upper index. The master function is $\Sigma_{\bs l}$-invariant.

A point $\bs t$ with complex coordinates is called {\it a critical
point} of $\Phi_{\bs d}(\,\cdot\,;\bs z)$ if the following system of
$l_1+\dots+l_{r}$ equations is satisfied
\bean\label{BAE}
-\sum_{s=1,\ s\neq j}^{l_1}
\frac{2}{t_j^{(1)}-t_s^{(1)}}+\sum_{s=1}^{l_2}
\frac{1}{t_j^{(1)}-t_s^{(2)}}=0\ , \notag
\\
-\sum_{s=1,\ s\neq j}^{l_i}\frac{2}{t_j^{(i)}-t_{s}^{(i)}}+\sum_{s=1}^{l_{i-1}}
\frac1{t_j^{(i)}-t_s^{(i-1)}}+\sum_{s=1}^{l_{i+1}}
\frac{1}{t_j^{(i)}-t_s^{(i+1)}}=0\ ,
\\
\sum_{s=1}^{n}
\frac 1 {t_j^{(r)}-z_s}
-\sum_{s=1,\ s\neq j}^{l_{r}}
\frac{2}{t_j^{(r)}-t_s^{(r)}}+\sum_{s=1}^{l_{r-1}}
\frac{1}{t_j^{(r)}-t_s^{(r-1)}}=0\ ,
\notag
\eean
where $j=1,\dots,l_1$ in the first group of equations,
$i=2,\dots,r-1$ and $j=1,\dots,l_i$ in the second group of
equations, $j=1,\dots,l_{r}$ in the last group of equations.

In other words, a point $\bs t$ is a critical point if
\bea
\left(\Phi_{\bs d}^{-1}
\frac{\partial \Phi_{\bs d} }{\partial t_j^{(i)}}\right)\!\!
(\bs t; \bs z)\ =\ 0\ ,
\qquad
i = 1 , \dots , r,\ j = 1 , \dots l_i\ .
\eea
In the Gaudin model, equations \Ref{BAE} are called {\it the Bethe
ansatz equations}.

The critical set is $\Sigma_{\bs l}$-invariant.

For a critical point $\bs t$, define the tuple $\bs y^{\bs t} = (y_1,\dots , y_{r})$
of polynomials in variable $x$,
\bean
\label{Tuple}
y_i(x) \ =\ \prod_{j=1}^{l_i}(x-t^{(i)}_j)\ ,
\qquad
i = 1, \dots , r\ .
\eean
Consider the linear differential operator of order $r+1$,
$$
D_{ \bs t} = ( \frac{d}{dx} -
\ln' ( \frac { T } { y_{r} } ) )
( \frac{d}{dx} - \ln' ( \frac {y_{r} } {y_{r-1} } ) )
\dots ( \frac{d}{dx} - \ln' ( \frac { y_2 }{ y_1 } ) )
( \frac{d}{dx} - \ln' ( y_1 ) ) ,
$$
where $\ln'(f)$ denotes $({df}/{dx})/f$ for any $f$. Denote by
$V_{\bs t}$ the kernel of $D_{\bs t}$.

Call $D_{\bs t}$ {\it the fundamental operator of the critical point $\bs t$}.
Call $V_{\bs t}$ {\it the fundamental space of the critical point $\bs t$}.

\begin{theorem}[Section 5 in \cite{MV2}]
\label{thm construction of spaces}
${}$

\begin{enumerate}
\item[$\bullet$]
The fundamental space $V_{\bs t}$ is an
$r+1$-dimensional space of polynomials with exponents $\bs d$ at infinity
and Wronskian $T$.
\item[$\bullet$]
The tuple $\bs y^{\bs t}$ can be recovered from the fundamental space
$V_{\bs t}$ as follows. Let
$f_1,\dots, f_{r+1}$ be a basis of $V_{\bs t}$, consisting of polynomials
with $\deg\, f_i = d_i$ for all $i$. Then $y_1,\dots,y_{r}$
are respective scalar multiples of the polynomials
$$
f_1 , \quad
{\rm Wr} \,(f_1,f_2) , \quad
{\rm Wr} \,(f_1,f_2,f_3) ,\quad
\dots\ ,
{\rm Wr}\, (f_1,\dots,f_{r})\ .
$$
\end{enumerate}
\end{theorem}

Thus distinct orbits of critical points define distinct $r+1$-dimensional
spaces $V$
with exponents $\bs d$ at infinity and Wronskian $T$.

\begin{theorem} [Theorem 6.1 in \cite{MV3}]
\label{theorem 4}
For generic complex $z_1,\dots,z_n$, the master function
$\Phi_{\bs d}(\,\cdot\,;\bs z)$
has $N(\bs d)$ distinct orbits of critical points.
\end{theorem}

Therefore, by Theorems \ref{thm construction of
spaces} and \ref{theorem 4}, we constructed $N(\bs d)$ distinct spaces of polynomials with
Wronskian $T$. All these spaces are fundamental spaces of critical
points of the master function $\Phi_{\bs d}(\,\cdot\,;\bs z)$.

\section{Bethe vectors}
\label{Bethe vectors}

\subsection{Generators} Let $E_{i,j}$, $i,j = 1, \dots, r+1$, be the standard generators
of $\glg_{r+1}$. The elements $E_{i,j}$, $i\neq j$, and $H_i =
E_{i,i}-E_{i+1,i+1}$, $i=1,\dots, r$, are the standard generators of
$\slg_{r+1}$. We have $\slg_{r+1} = \n_+\oplus\h\oplus\n_-$ where
\bea &
\n_+ = \oplus_{i<j} \C\cdot E_{i,j}\ , \qquad \h = \oplus_{i=1}^r
\C\cdot H_i\ , \qquad \n_- = \oplus_{i>j} \C\cdot E_{i,j}\ .
\eea

\subsection{Construction of Bethe vectors}
\label{construction of B vectors}
For $\mu \in \h^*$, denote by $L_{\om_r}^{\otimes n}[\mu]$ the
vector subspace of $L_{\om_r}^{\otimes n}$ of vectors of weight
$\mu$ and by $\sing\,L_{\om_r}^{\otimes n}[\mu]$
the vector subspace of
singular vectors of weight $\mu$,
\bea
L_{\om_r}^{\otimes n}[\mu] &= & \{
v \in L_{\om_r}^{\otimes n}\ {} |\ {} hv\, = \,\langle \mu, h \rangle \,v
\ {\rm for\ any }\ h\in \h \}\
\\
\sing\,L_{\om_r}^{\otimes n}[\mu] & = & \{
v \in L_{\om_r}^{\otimes n}\ {} |\ {} \n_+v\,=\,0, \ {} hv\, = \,\langle \mu, h \rangle \,v
\ {\rm for\ any }\ h\in \h \}\ .
\eea
For a given $\bs l = (l_1,\dots,l_r)$, set $l = l_1 + \dots + l_r$ and
\bea
\mu\ =\ n\om_r-\sum_{i=1}^rl_i\al_i\ .
\eea
Let $\C^l$ be the space with coordinates $t^{(i)}_j$, $i=1,\dots,r$, $j=1,\dots, l_i$,
and $\C^n$ the space with coordinates $z_1,\dots, z_n$.
We construct a rational map
\bea
\omega \ :\ \C^l \times \C^n\ \to \nL
\eea
called {\it the universal weight function}.

Let $P(\bs l,n)$ be the set of sequences $I\ = \ (i_1^1, \dots ,
i^1_{k_1};\ \dots ;\ i^n_1, \dots , i^n_{k_n})$ of integers in $\{1,
\dots , r\}$ such that for all $i = 1, \dots , r$, the integer $i$
appears in $I$ precisely $l_i$ times. For $I\in P(\bs l,n)$, the $l$
positions in $I$ are partitioned into subsets $I_1,\dots,I_r$, where
$I_i$ consists of positions of the integer $i$. Fix a labeling of
positions in $I_i$ by numbers $1,\dots,l_i$. Thus to every position
${}^a_b$ in $I$ we assign an integer, the labeling number of this
position in the corresponding subset. Denote this integer by
$j(I^a_b)$. For $\sigma = (\sigma_1,\dots,\sigma_r)\in \Sigma_{\bs
l}$, denote by $t(I^a_b;\sigma)$ the variable
$t^{(i)}_{\sigma_{i}(j)}$ where $i=i^a_b$, $j=j(I^a_b)$, and
$\sigma_{i}(j)$ denotes the image of $j$ under the permutation
$\sigma_i$. For a given $\sigma$, the assignment of this variable to
a position establishes a bijection of $l$ positions of $I$ and
the set $\{t^{(1)}_{1},$ $\dots,$ $t_{l_1}^{(1)},$ $\dots,
t^{(r)}_{1},\dots,t_{l_{r}}^{(r)}\}$.

Fix a highest weight vector $v_{\omega_r}$ in $L_{\om_r}$. To every $I \in P(\bs
l, n)$, assign the vector
\bea
E_I v\ =\ E_{i_1^1+1\,,\,i_1^1} \dots
E_{i_{k_1}^1+1\,,\,i_{k_1}^1}v_{\omega_r} \otimes \dots \otimes E_{i_1^n+1\,,\,i_1^n}
\dots E_{i_{k_n}^n+1\,,\,i_{k_n}^n}v_{\omega_r} \
\eea
in $\nL$ and scalar functions
$\omega_{I,\sigma}$ labeled by $\sigma = (\sigma_1,\dots,\sigma_r)\in
\Sigma_{\bs l}$, where
\bea
\omega_{I,\sigma} \ =\
\omega_{I,\sigma,1}(z_1) \dots \omega_{I,\sigma,n}(z_n)
\eea
and
\bea
\omega_{I,\sigma,j}(z_j) \ =\ \frac 1 {(t(I^j_1,\sigma) -
t(I^j_2,\sigma))} \ \dots\ \frac 1 {(t(I^j_{k_j-1},\sigma) -
t(I^j_{k_j},\sigma))} \ \frac 1 {(t(I^j_{k_j},\sigma) - z_j)} \ .
\eea
We set
\bean
\label{bethe vector}
\omega(\bs t; \bs z) \ =\
\sum_{I\in P(\bs l,n)}\ \sum_{\sigma\in \Sigma_{\bs l}}\
\omega_{I,\sigma}\ E_I v\ .
\eean The universal weight function is
symmetric with respect to the $\Sigma_{\bs l}$-action on variables
$t^{(i)}_j$.

\medskip
Other formulas for the universal weight function see in
\cite{M, RSV}.

The universal weight function was introduced in \cite{SV} to solve the
KZ equations.

\medskip

\noindent
{\bf Examples.} If $n=2$ and $\bs l = (1, 1, 0, \dots , 0)$, then
\bea
\omega(\bs t;\bs z) =
\frac {1} {(t^{(1)}_1-t^{(2)}_1)
( t^{(2)}_1
- z_1)}
E_{2,1} E_{3,2} v_{\om_r} \otimes v_{\om_r} +
\frac 1{(t^{(2)}_1-t^{(1)}_1)(t^{(1)}_1 -
z_1)}E_{3,2}E_{2,1}v_{\om_r}\otimes v_{\om_r}\kern-2em
\\
+
\frac 1{(t^{(1)}_1-z_1)(t^{(2)}_1-z_2)}E_{2,1}v_{\om_r}\otimes E_{3,2}v_{\om_r} +
\frac 1{(t^{(2)}_1-z_1)(t^{(1)}_1-z_2)}E_{3,2}v_{\om_r}\otimes E_{2,1}v_{\om_r}
\kern-2em
\\
+
\frac 1{(t^{(1)}_1-t^{(2)}_1)(t^{(2)}_1-z_2)}v_{\om_r}\otimes E_{2,1}E_{3,2}v_{\om_r} +
\frac 1{(t^{(2)}_1-t^{(1)}_1)(t^{(1)}_1-z_2)}v_{\om_r}\otimes E_{3,2}E_{2,1}v_{\om_r}\ .\kern-2.6em
\eea
If $\bs l = (2, 0, \dots , 0)$, then
\bea
\omega(\bs t;\bs z)
&=&
(\frac 1{(t_1^{(1)}-t_2^{(1)})(t_2^{(1)}-z_1)}+\frac 1{(t_2^{(1)}-t_1^{(1)})(t_1^{(1)}-z_1)})
\ {}E_{2,1}^2v_{\om_r}\otimes v_{\om_r}
\\
& +&
(\frac 1{(t_1^{(1)}-z_1)(t_2^{(1)}-z_2)}+\frac 1{(t_2^{(1)}-z_1)(t_1^{(1)}-z_2)})
\ {} E_{2,1}v_{\om_r}\otimes E_{2,1}v_{\om_r}
\\
& +&
(\frac 1{(t_1^{(1)}-t_2^{(1)})(t_2^{(1)}-z_2)}+\frac 1{(t_2^{(1)}-t_1^{(1)})(t_1^{(1)}-z_2)})
\ {} \ {}\ v_{\om_r}\otimes E_{2,1}^2v_{\om_r}\ .
\eea

The values of the universal weight function at the critical points of
the master function are called {\it the Bethe vectors}.

The Bethe vectors of critical points of the same $\Sigma_{\bs l}$-orbit
coincide, since both the critical point equations and the universal
weight functions are $\Sigma_{\bs l}$-invariant.

The universal weight function takes values in $\nL$. But if $\bs t$ is
a critical point of the master function, then the Bethe vector
$\om(\bs t; \bs z)$ belongs to the subspace of singular vectors
$\snL\subset\nL$, see \cite{RV}, cf. comments on this fact in Section
2 of \cite{MV3}.

\medskip

By Theorem \ref{theorem 4}, the master function $\Phi_{\bs d}(\,\cdot\,;\bs z)$ has
$N(\bs d)$ distinct orbits of critical points for generic $\bs z$.
Choose a representative in each of the orbits: \ $\bs t^1, \dots , \bs
t^{N(\bs d)}$. These critical points define the collection of Bethe
vectors:\ $\om(\bs t^1; \bs z), \dots , \om(\bs t^{N(\bs d)}; \bs z) \in
\snL$. The space $\snL$ has dimension $N(\bs d)$.

\begin{theorem}
[Theorem 6.1 in \cite{MV3}]
\label{theorem from MV3 about basis}
For generic $\bs z$, the Bethe vectors form a basis in $\snL$.
\end{theorem}

\subsection{The Gaudin model} The Gaudin Hamiltonians on $\snL$ is
a collection of linear operators acting on $\snL$ and (rationally) depending
on a complex parameter $x$. We use the construction of the Gaudin Hamiltonians
suggested in \cite{T, CT}, see also \cite{MTV}.
We consider the $\slg_{r+1}$-module $L_{\om_r}$
as the $\glg_{r+1}$-module of highest weight $(0,\dots,0,-1)$.

To define the Gaudin Hamiltonians, for all
$i,j = 1, \dots, r+1$, consider the differential operators
$$
X_{i,j}(x)\ =\ \delta_{i,j}\,\frac{d}{dx}
\ -\,\sum_{s=1}^n\,\frac{E_{j,i}^{(s)}}{x-z_s}
$$
where $\delta_{i,j}$ is the Kronecker symbol and
$E_{j,i}^{(s)}=\,1^{\otimes(s-1)}\,\otimes \,E_{j ,i}\otimes\,
1^{\otimes(n-s)}$.
These differential operators act on
$L_{\om_r}^{\otimes n}$-valued functions in $x$. The order of $X_{ij}$ is one
if $i=j$ and is zero otherwise.

Set
\beq
\label{Mbold}
\bs M\,\,=\sum_{\sigma \in \Sigma_{r+1}}\ (-1)^{\sigma}\
X_{1 , \sigma(1)}(x)\ X_{2 , \sigma(2)}(x)\ \dots\ X_{r+1 , \sigma(r+1)}(x)
\eeq
where
$(-1)^{\sigma}$ denotes the sign of the permutation.

For example, for $r=1$, we have
\bea
\bs M \ =\ (\, \frac{d}{dx}
\ - \,\sum_{s=1}^n\,\frac{E_{1,1}^{(s)}}{x-z_s}\, )\
(\, \frac{d}{dx}
\ - \,\sum_{s=1}^n\,\frac{E_{2,2}^{(s)}}{x-z_s} \,)\
-\
( \,\sum_{s=1}^n\,\frac{E_{2,1}^{(s)}}{x-z_s} \,)\
( \,\sum_{s=1}^n\,\frac{E_{1,2}^{(s)}}{x-z_s} \,)\ .
\eea

Write
$$
\bs M\,=\,\frac{d^{r+1}}{dx^{r+1}}+M_1(x)\frac{d^r}{dx^r}+\dots+M_{r+1}(x)\ ,
$$
where $M_i(x) : \nnL \to \nnL$ are linear operators depending on $x$.
The coefficients $M_1(x)$, \dots, $M_{r+1}(x)$ are called {\it the Gaudin Hamiltonians}.

It is the well-known that
\begin{enumerate}
\item[$\bullet$]
The Gaudin Hamiltonians commute: $[M_i(u),M_j(v)]=0$ for all $i,j,u,v$,

\item[$\bullet$]
The Gaudin Hamiltonians commute with the $\glg_{r+1}$-action on $\nnL$,
\newline
in particular, they preserve $\snL$.
\end{enumerate}
The first statement see, for example, in
\cite{KuS, T, CT}, Proposition 7.2 \cite{MTV}.
The second statement see, for example, in \cite{KuS} and
Proposition 8.3 in \cite{MTV}.

\begin{theorem}
[Theorem 9.2 in \cite{MTV}]
\label{eigenvalue theorem}
For any critical point $\bs t$ of the master function
$\Phi_{\bs d}(\,\cdot\,;\bs z)$, the Bethe vector
$\btz$ is an eigenvector of $M_1(x),\dots, M_{r+1}(x)$. The corresponding
eigenvalues $\mu_1(x),\dots,\mu_{r+1}(x)$ are given by the formula
\bea
&& \frac{d^{r+1}}{dx^{r+1}}+\mu_1(x)\frac{d^r}{dx^r}+\dots+\mu_{r+1}(x)\,={}
\\
&&
\phantom{aaa}\,(\frac{d}{dx}+\ln'(y_1))\;(\frac{d}{dx}+\ln'(\frac{y_2}{y_1}))\;\dots
(\frac{d}{dx}+\ln'(\frac{y_r}{y_{r-1}}))\;(\frac{d}{dx}+\ln'(\frac{T}{y_r}))\ .
\eea
\end{theorem}

Set
\bea
\bs K\kern-16pt
&& {}=\;
\frac{d^{r+1}}{dx^{r+1}}\ -\ \frac{d^r}{dx^r} M_1(x)\ +
\ \dots + (-1)^{r+1}M_{r+1}(x)
\\
&&{}=\;
\frac{d^{r+1}}{dx^{r+1}}\ +\ K_1(x)\frac{d^r}{dx^r}\ +\ \dots+K_{r+1}(x) \ .
\eea
This is the differential operator that is formally conjugate to
the differential operator
$(-1)^{r+1}\bs M$. The coefficients
$K_i(x) : \nnL \to \nnL$ are linear operators depending on $x$.
These coefficients can
be expressed as differential polynomials in $M_1(x),$ $\dots,$ $M_{r+1}(x)$.
For instance,
\bea
K_1(x)\,=\,-\,M_1(x)\,,
\qquad
K_2(x)\,=\,M_2(u)\,-\,r\,\frac{d}{dx} M_1(x)\,,
\eea
and so on. The operators $K_1(x), \dots, K_{r+1}(x)$
pairwise commute, \
\ $[K_i(u),K_j(v)]=0$ for all $i,j,u,v$,
\ and commute with the $\glg_{r+1}$-action on $\nnL$.

The operators $K_1(x), \dots, K_{r+1}(x)$ are also called {\it the Gaudin Hamiltonians}.

\medskip

For any critical point $\bs t$ of the master function
$\Phi_{\bs d}(\,\cdot\,;\bs z)$, the Bethe vector $\btz$ is an eigenvector of
the Gaudin Hamiltonians
$K_1(x),\dots, K_{r+1}(x)$. The corresponding eigenvalues
$\la_1(x),\dots,\la_{r+1}(x)$ are given by the formula
\bea
&& \frac{d^{r+1}}{du^{r+1}}+\la_1(x)\frac{d^r}{dx^r}+\dots+\la_{r+1}(x)\,
\\
&&
\phantom{aa}=\,(\frac{d}{dx}-\ln'(\frac{T}{y_r}))\;
(\frac{d}{dx}-\ln'(\frac{y_r}{y_{r-1}}))\;\dots
(\frac{d}{dx}-\ln'(\frac{y_2}{y_1}))\;(\frac{d}{dx}-\ln'(y_1))\ .
\eea
Notice that this is the fundametal differential operator
$D_{\bs t}$ of the critical point $\bs t$.

\begin{cor}
\label{cor on simple spectrum}
For generic $\bs z$,
\begin{enumerate}
\item[$\bullet$]
the Bethe vectors form an eigenbasis of the
Gaudin Hamiltonians $K_1(x)$, \dots, $K_{r+1}(x)$,
\item[$\bullet$]
the Gaudin Hamiltonians $K_1(x), \dots, K_{r+1}(x)$
have simple spectrum, that
is the eigenvalues of the Gaudin Hamiltonians separate the basis Bethe
eigenvectors.
\end{enumerate}
\end{cor}

The first statement of the corollary follows from Theorem
\ref{theorem from MV3 about basis} and Theorem \ref{eigenvalue theorem}.

The second statement of the corollary follows from the fact that if
two Bethe vectors had the same eigenvalues, then they would have the
same fundamental operators, hence the same fundamental spaces. But
the fundamental space uniquely determines the orbit of the
corresponding critical point by Theorem \ref{thm construction of
spaces}. Hence the two Bethe vectors correspond to the same orbit
of critical points. Hence the two Bethe vectors are equal.

\subsection{The Shapovalov form and real $\bs z$}
\label{sec The Shapovalov }
Define the anti-involution $\tau : \glg_{r+1} \to \glg_{r+1} $ by sending
$E_{i,j}$ to $E_{j,i}$ for all $i,j$.

Let $W$ be a highest weight $\glg_{r+1}$-module with highest weight vector $w$.
{\it The Shapovalov form} on $W$ is the unique
symmetric bilinear form $S$
defined by the conditions:
\bea
S(w, w) = 1 ,
\qquad
S(gu, v) = S(u, \tau(g)v)
\eea
for all $u,v \in W$ and $g \in \glg_{r+1}$, see \cite{K}.
The Shapovalov form is non-degenerate on an irreducible module
$W$ and is positive definite on the real part of the irreducible module
$W$.

Let $L_{\La_1}\!\otimes \dots \otimes L_{\La_n}$ be the tensor product
of irreducible highest weight $\glg_{r+1}$-modules. Let $v_{\La_i} \in L_{\La_i}$
be a highest weight vector and $S_i$ the corresponding Shapovalov form
on $L_{\La_i}$. Define the symmetric bilinear form on the tensor
product by the formula $S = S_1 \otimes \dots \otimes S_n.$ The form
$S$ is called {\it the tensor Shapovalov form}.

\begin{theorem}[Proposition 9.1 in \cite{MTV}]
\label{shap and gaudin}
The Gaudin hamiltonians $K_1(x), \dots , K_{r+1}(x)$ are symmetric
with respect to the tensor Shapovalov form $S$,
$$
S(K_i(x)u, v)\ =\ S(u, K_i(x)v)
\qquad
{\rm for\ all} \quad i, x, u, v\ .
$$
\end{theorem}

\begin{cor}
If all of $z_1,\dots,z_n, x$ are real numbers, then the Gaudin Hamiltonians
$K_1(x)$, \dots , $K_{r+1}(x)$ are real linear operators on the real part of the
tensor product $L_{\La_1}\!\otimes \dots \otimes L_{\La_n}$. These operators
are symmetric with respect to the positive definite tensor Shapovalov form. Hence
they are simultaneously diagonalizable and have real spectrum.
\end{cor}

\subsection{Proof of Theorem \ref{theorem 3}}
If $z_1,\dots,z_n, x$ are real, then the Gaudin Hamiltonians on $\snL$
have real
spectrum as symmetric operators on a Euclidean space. If $\bs t$ is a
critical point of $\Phi_{\bs d}(\,\cdot\,;\bs z)$
and $\la_1(x), \dots, \la_{r+1}(x)$ are the eigenvalues
of the corresponding Bethe vector $\btz$, then $\la_1(x), \dots,
\la_{r+1}(x)$ are real rational functions. Hence the fundamental
differential operator $D_{\bs t}$ has real coefficients. Therefore,
the fundamental vector space of polynomials $V_{\bs t}$ is real. Thus for generic
real $z_1,\dots,z_n$ we have $N(\bs d)$ distinct real spaces of polynomials
with exponents $\bs d$ at infinity and Wronskian
$\prod_{s=1}^n(x-z_s)$. Theorem \ref{theorem 3} is proved.

\section{Appendix A}
\subsection{The differential operator $\bs K$ has polynomial solutions only}
Let $z_1,\dots,z_n\in\C$. Let $\La_1,\dots,\La_n,\La_\infty\in\h^*$ be dominant
integral weights. Assume that the irreducible $\slg_{r+1}$-module
$L_{\La_\infty}$ is a submodule of the tensor product
$L_{\La_1}\!\otimes\dots\otimes L_{\La_n}$.

For any $s=1,\dots,n,\infty$, and $i= 1,\dots, r$, set
$m_{s,i}=(\La_s\,,\sum_{j=1}^i\al_j)$ and
\bea
l\;=\;\frac1{r+1}\;\sum_{i=1}^r\,\Bigl(\;\sum_{s=1}^n\, m_{s,i}\ -
\ m_{\infty,i}\Bigr)\ .
\eea

For any $s=1,\dots,n$, we will consider the $\slg_{r+1}$-module $L_{\La_s}$
as the $\glg_{r+1}$-module of highest weight
$\bigl(0,-m_{s,1},-m_{s,2},\dots,-m_{s,r}\bigr)$.
The $\slg_{r+1}$-module $L_{\La_\infty}$, as a submodule of
the $\glg_{r+1}$-module $L_{\La_1}\!\otimes\dots\otimes L_{\La_n}$, has
the $\glg_{r+1}$-highest weight
$$
(-l,-l-m_{\infty,1},-l-m_{\infty,2},\dots,-l-m_{\infty,r})\;.
$$

\begin{theorem}
\label{monodromy}
{\strut}
\begin{enumerate}
\item[(i)]
Consider the operator $\bs K$ as a differential operator acting on
$L_{\La_1}\!\otimes\dots\otimes L_{\La_n}\!$-valued functions in $x$.
Then all singular points of the operator $\bs K$ are regular and lie in the set
$\{z_1, \dots , z_n,$ $\infty\}$.
\item[(ii)]
Let $u(x)$ be any germ of an
$L_{\La_1}\!\otimes\dots\otimes L_{\La_n}\!$-valued function
such that $\bs Ku = 0$. Then $u$ is the germ of an
$L_{\La_1}\!\otimes\dots\otimes L_{\La_n}\!$-valued polynomial in $x$.
\item[(iii)]
Let $w \in \sing\,(L_{\La_1}\!\otimes\dots\otimes L_{\La_n})[\La_\infty]$ be an
eigenvector of the operators $K_1(x),\dots,$ $K_{r+1}(x)$ with the eigenvalues
$\la_1(x),\dots,\la_{r+1}(x)$, respectively. Consider the scalar
differential operator
$$
D_w\;=\;\frac{d^{r+1}}{dx^{r+1}}\ +\ \la_1(x)\;\frac{d^r}{dx^r} +
\ \dots+\la_{r+1}(x)\;.
$$
Then the exponents of the differential operator $D_w$ at $\infty$ are
$\;-l\,,\,-m_{\infty,1}-1-l\,,\,\dots\,,$ $-m_{\infty,r}-r-l$.
\item[(iv)]
If $z_1,\dots,z_n$ are distinct, then
for any ${s=1,\dots,n}$, the exponents of the differential operator
$D_w$ at $z_s$ are $\;0\,,\allowbreak\,m_{s,1}+1\,,\,\dots\,,\,m_{s,r}+r$.
\item[(v)]
The kernel of the differential
operator $D_w$ is an $r+1$-dimensional space of polynomials.
\end{enumerate}
\end{theorem}

\begin{proof}
Part (i) is a direct corollary of the definition of the operator $\bs K$.

We first prove part (ii) in the special
case of $\La_1=\dots=\La_n=\om_r$ and generic $z_1,\dots, z_n$.
By construction, the operator $\bs K$ commutes with the
$\glg_{r+1}$-action on $L_{\om_r}^{\otimes n}$.
This fact and Theorem~\ref{theorem from MV3
about basis} imply that $\bs K$ has an eigenbasis consisting of the
Bethe vectors and their images under the $\glg_{r+1}$-action. Then
by Theorems~\ref{eigenvalue theorem} and \ref{thm construction of
spaces}, all solutions of the differential equation $\bs Ku=0$ are
polynomials.

The proof of part (ii) for arbitrary
$\La_1, \dots , \La_n$ and $z_1,\dots,z_n$ clearly follows from the special case and
the following remarks:
\begin{enumerate}
\item[$\bullet$]
The operator $\bs K$ is well-defined for any $z_1,\dots,z_n$,
not necessarily distinct, and rationally depends on $z_1,\dots,z_n$.
\item[$\bullet$]
If for generic $z_1,\dots,z_n$, all solutions of the differential
equation $\bs Ku=0$ are polynomials, then for any $z_1,\dots,z_n$, all
solutions of the differential equation $\bs Ku=0$ are polynomials.
\item[$\bullet$] Assume that some of $z_1,\dots,z_n$
coincide. Partition the set $\{z_1,\dots,z_n\}$ into several groups of
coinciding points of sizes $n_1, \dots, n_k$,
$n_1+\dots+n_k=n$. Denote the representatives in the groups by
$u_1,\dots,u_k \in \C$ where $u_1,\dots,u_k$ are distinct.
Denote $W_s=L_{\om_r}^{\otimes n_s}$ for
$s=1,\dots,k$. Choose an irreducible module $L_{\nu_s} \subset W_s$
for every $s$. Then the operator $\bs K$ defined for those
$z_1,\dots,z_n$ on $W_1\otimes \dots\otimes W_k$ preserves the
space of functions with values in the
submodule $L_{\nu_1}\otimes \dots \otimes L_{\nu_k}$. If we
restrict $\bs K$ to the space of functions with values in
$L_{\nu_1}\otimes \dots \otimes L_{\nu_k}$, then this restriction
coincides with the operator $\bs K$ defined for the tensor product
$L_{\nu_1}\otimes \dots \otimes L_{\nu_k}$ and $u_1,\dots,u_k$.
\item[$\bullet$]
Any highest weight irreducible module
is a submodule of a suitable tensor power of $L_{\om_r}$.
\end{enumerate}
Part (ii) is proved.

In order
to calculate the exponents of the operator $D_w$ at singular points, we calculate the
exponents of its formal conjugate operator. Namely, we consider the
operator
\bea
D_w^*\kern-18pt &&
{}=\;\frac{d^{r+1}}{dx^{r+1}}\ -\ \frac{d^r}{dx^r}\;\la_1(x)\ +
\ \dots+(-1)^{r+1}\la_{r+1}(x)
\\
&& {}=\;\frac{d^{r+1}}{dx^{r+1}}\ +\ \mu_1(x)\;\frac{d^r}{dx^r} +
\ \dots+\mu_{r+1}(x)\ .
\eea
The vector $w$ is an eigenvector of the operators
$M_1(x),\dots,M_{r+1}(x)$ with the eigenvalues
$\mu_1(x),\dots, \mu_{r+1}(x)$, respectively.

\begin{lem}
\label{exponents}
Let the exponents of $D_w^*$ at a point $x=z$ be $p_1,\dots,p_{r+1}$.
Then the exponents of $D_w$ at the point $x=z$ are $r-p_{r+1},\dots,r-p_1$.
\hfill
$\square$
\end{lem}

Consider the following $U(\glg_{r+1})$-valued polynomial
\bean
\label{Ax}
&&
A(x)\ =\ \sum_{\sigma\in \Sigma_{r+1}} (-1)^{\sigma}\
\bigl((x-r)\,\delta_{1,\,\sigma(1)}-E_{\sigma(1),1}\bigr)\ \dots
\phantom{aaaaaaaaa}
\\
&&
\phantom{aaaaaaaaaaaaa}
\dots\
\bigl((x-1)\,\delta_{r,\,\sigma(r)}-E_{\sigma(r),r}\bigr)
\bigl(x\,\delta_{r+1,\,\sigma(r+1)}-E_{\sigma(r+1),r+1}\bigr)\ .
\notag
\eean
It is known that the coefficients of this polynomial are central elements
in $U(\glg_{r+1})$, see, for example, Remark~2.11 in \cite{MNO}.
If $v$ is a singular vector of a $\glg_{r+1}$-weight $(p_1,\dots, p_{r+1})$,
then formula \Ref{Ax} yields
\bea
A(x)\,v\ =\ \prod_{i=1}^{r+1}\,(x-r-1+i-p_i)\,v\ .
\eea
Hence, the operator $A(x)$ acts on $L_{\La_s}$ as the identity operator
multiplied by
\bea
\psi_s(x)\;=\;\prod_{i=0}^r\,(x-r+i+m_{s,i})\;.
\eea

\medskip
Let $s=1,\dots,n$.
It follows from formula \Ref{Mbold} that the indicial polynomial of $D_w^*$
at the singular point $z_s$ is the eigenvalue of the operator
$1^{\otimes(s-1)}\otimes A(x)\otimes 1^{\otimes(n-s)}$ on the vector $w$,
that is, $\psi_s(x)$. Similarly, the indicial polynomial of $D_w^*$ at infinity
is the eigenvalue of $A(-x)$ acting on the vector $w$ which belongs to
the submodule $L_{\La_\infty}$ of the $\glg_{r+1}$-module
$L_{\La_1}\!\otimes\dots\otimes L_{\La_n}$ , that is,
$$
\psi_\infty(x)\;=\;\prod_{i=0}^r\,(-x-r+i+l+m_{\infty,i})\;.
$$
Hence, by Lemma~\ref{exponents}, the exponents of the operator $D_w$ are
as required. This proves parts (iii) and (iv). Part (v) follows from parts (i-iv).
\end{proof}

\begin{cor}
\label{last cor}
Assume that the operators $K_1(x),\dots,K_{r+1}(x)$ acting on
the subspace of weight singular vectors
$\sing(L_{\La_1}\!\otimes\dots\otimes L_{\La_n})[\La_\infty]$ are
diagonalizable and have simple joint spectrum. Then there exist
$\dim \,\sing(L_{\La_1}\!\otimes\dots\otimes
L_{\La_n})[\La_\infty]$ distinct polynomial $r+1$-dimensional
spaces $V$ with the
following properties.
If $D$ is the fundamental differential operator of such a space, then
$D$ has singular points at $z_1,\dots,z_n,\infty$ only, with the exponents
$\,0\,,\,m_{s,1}+1\,,\,\dots\,,\,m_{s,r}+r$ at $z_s$ for any $s$, and
the exponents
$\,-l\,,\,-m_{\infty,1}-1-l\,,\,\dots\,,\,-m_{\infty,r}-r-l$ at $\infty$.
\end{cor}

Consider all $r+1$-dimensional
polynomial spaces $V$, whose fundamental operator has exponents at
$z_1,\dots,z_n,\infty$ indicated in Corollary \ref{last cor}. Schubert calculus
says that the number of such spaces is not
greater than the dimension of $\sing(L_{\La_1}\!\otimes\dots\otimes
L_{\La_n})[\La_\infty]$, see for example \cite{MV2}. Thus, according
to Corollary \ref{last cor}, the simplicity of the spectrum of the
Gaudin Hamiltonians on $\sing\,(L_{\La_1}\!\otimes\dots\otimes
L_{\La_n})[\La_\infty]$ implies the transversality of the Schubert
cycles corresponding to these exponents at $z_1,\dots, z_n, \infty$,
cf. \cite{MV2} and \cite{EH}.
\medskip

Recall also that the operators $K_1(x),\dots,K_{r+1}(x)$ acting on
$\sing(L_{\La_1}\!\otimes\dots\otimes L_{\La_n})[\La_\infty]$ are
diagonalizable if $z_1,\dots,z_n$ are real,
see Section \ref{sec The Shapovalov }.

\medskip

\begin{rem}
It was conjectured in \cite{CT} that the monodromy of
the differential operator $\bs M$, acting on
$L_{\La_1}\!\otimes\dots\otimes L_{\La_n}$-valued functions in $x$, is
trivial. However, the proof of this statement in \cite{CT} is not
satisfactory. On the other hand
Theorem~\ref{monodromy} implies that the monodromy of
the differential operator $\bs K$, acting on
$L_{\La_1}\!\otimes\dots\otimes L_{\La_n}$-valued functions in $x$, is
trivial. Together with
Theorem~\ref{shap and gaudin}, this implies that the monodromy of the
operator $\bs M$ is trivial as well.
\end{rem}

\subsection{Bethe vectors in
$\sing\,(L_{\La_1}\!\otimes\dots\otimes L_{\La_n})[\La_\infty]$}
Let $\bs z=(z_1,\dots,z_n)\in \C^n$ be a point with distinct coordinates.
Let $\La_1,\dots,\La_n,\La_\infty\in\h^*$ be dominant
integral weights. Assume that the irreducible $\slg_{r+1}$-module
$L_{\La_\infty}$ is a submodule of the tensor product
$L_{\La_1}\!\otimes\dots\otimes L_{\La_n}$.

Introduce $\bs l = (l_1,\dots,l_r)$ by the formula $\La_\infty =
\sum_{s=1}^n\La_s-\sum_{i+1}^rl_i\al_i$. Set $l = l_1+\dots+l_r$. Consider
the associated master function
\beq
\label{sl master}
\Phi(\bs t;\bs z)\ =\ \prod_{i=1}^r \prod_{j=1}^{l_i}\prod_{s=1}^n
(t_j^{(i)}-z_s)^{-(\La_s,\al_i)} \prod_{i=1}^r\prod_{1\leq j<s\leq
l_i} (t_j^{(i)}-t_s^{(i)})^{2}
\prod_{i=1}^{r-1}\prod_{j=1}^{l_i}\prod_{k=1}^{l_{i+1}}
(t_j^{(i)}-t_k^{(j+1)})^{-1} \ . \notag
\eeq
Consider the universal
weight function $\om : \C^l\times\C^n \to
(L_{\La_1}\!\otimes\dots\otimes L_{\La_n})[\La_\infty]$ defined by
the formulas of Section \ref{construction of B vectors}. The value
$\om(\bs t; \bs z)$ of the universal weight function at a critical
point $\bs t$ of the master function $\Phi(\,\cdot\,;\bs z)$ is called {\it a Bethe vector},
see \cite{RV, MV2}. The Bethe vector belongs to
$\sing\,(L_{\La_1}\!\otimes\dots\otimes L_{\La_n})[\La_\infty]$, see
\cite {RV}.

For a critical point $\bs t$, define the tuple $\bs y^{\bs t} = (y_1,\dots , y_{r})$
of polynomials in variable $x$ by formulas of Section \ref{construction of spaces section}.
Define polynomials $T_1,\dots,T_r$ in $x$ by the formula
\bea
&
T_i(x) \ =\ \prod_{s=1}^n (x-z_s)^{(\La_s,\al_i)}\ .
\eea
Consider the linear differential operator of order $r+1$,
$$
D_{ \bs t} = ( \frac{d}{dx} -
\ln' ( \frac { T_1\dots T_r } { y_{r} } ) )\
( \frac{d}{dx} - \ln' ( \frac {y_{r}T_1\dots T_{r-1} } {y_{r-1} } ) )
\dots ( \frac{d}{dx} - \ln' ( \frac {y_2 T_1}{ y_1 } ) ) \
( \frac{d}{dx} - \ln' ( y_1 ) ) \ .
$$
All singular points of $D_{\bs t}$ are regular and lie in
$\{ z_1, \dots , z_n,\infty \}$. The exponents of $D_{\bs t}$ at $z_s$
are $\,0\,,\,m_{s,1}+1\,,\,\dots\,,\,m_{s,r}+r$ for any $s$, and
the exponents of $D_{\bs t}$ at $\infty$ are
$\,-l\,,\,-m_{\infty,1}-1-l\,,\,\dots\,,\,-m_{\infty,r}-r-l$.
The kernel $V_{\bs t}$ of $D_{\bs t}$ is an
$r+1$-dimensional space of polynomials, see \cite{MV2}.

The tuple $\bs y^{\bs t}$ can be recovered from $V_{\bs t}$ as follows. Let
$f_1,\dots, f_{r+1}$ be a basis of $V_{\bs t}$, consisting of monic polynomials
of strictly increasing degree. Then $y_1,\dots,y_{r}$
are respective scalar multiples of the polynomials
$$
f_1\ ,\quad
\frac{{\rm Wr} \,(f_1,f_2)}{T_1}\ , \quad
\frac{{\rm Wr} \,(f_1,f_2,f_3)}{T_2T_1^2}\ ,
\quad\dots\quad ,
\ \frac{{\rm Wr}\, (f_1,\dots,f_{r})}{T_{r-1}T_{r-2}^2\dots T_1^{r-1}}\ ,
$$
see \cite{MV2}.
\begin{theorem}[Theorem 8.2 in \cite{MTV}]
\label{THM}
For any
critical point $\bs t$ of the master function $\Phi(\,\cdot\,;\bs z)$, the Bethe vector
$\btz$ is an eigenvector of $K_1(x),\dots, K_{r+1}(x)$ and the corresponding
eigenvalues $\la_1(x),\dots,\la_{r+1}(x)$ are given by the formula
\bea
\frac{d^{r+1}}{du^{r+1}}+\la_1(x)\frac{d^r}{dx^r}+\dots+\la_{r+1}(x) =
D_{\bs t}\ .
\eea
\end{theorem}

\begin{cor}
Any two distinct nonzero
Bethe vectors cannot have the same eigenvalues for all Gaudin Hamiltonians.
\end{cor}

The proof of the corollary is similar to the proof of the second statement of
Corollary \ref{cor on simple spectrum}.

\section{Appendix B}

Let $\g$ be a simple Lie algebra, $\h$ its Cartan
subalgebra, $\al_1, \dots, \al_r \in \h^*$ simple roots, $(\,,\,)$ the standard
invariant scalar product on $\g$. Let $\bs \La = (\La_1, \dots ,
\La_n)$ be integral dominant weights of $\g$. Let $\bs l =
(l_1,\dots,l_r)$ be non-negative integers such that the
weight
$\La_\infty = \sum_{s=1}^n\La_s - \sum_{i=1}^rl_i\al_i$ is dominant
integral. Let $\bs z = (z_1,\dots,z_n)$ be distinct complex numbers.
Introduce the associated {\it master function} of variables
\ $
\bs t = (t^{(1)}_{1},\dots,t_{l_1}^{(1)},\dots,
t^{(r)}_{1},\dots,t_{l_r}^{(r)})\
$
by the formula
\bea
\lefteqn{\Phi_{\g,\bs \La,\bs l} (\bs t;\bs z) = }
\notag
\\
&& \prod_{i=1}^r\prod_{j=1}^{l_i}\prod_{s=1}^n
(t_j^{(i)}-z_s)^{-(\La_s,\al_i)}
\prod_{i=1}^r\prod_{1\leq j<s\leq l_i} (t_j^{(i)}-t_s^{(i)})^{(\al_i,\al_i)}
\!\!\!
\prod_{1\leq i<j\leq r}\prod_{s=1}^{l_i}\prod_{k=1}^{l_{j}}
(t_s^{(i)}-t_k^{(j)})^{(\al_i,\al_j)} \ .
\eea
The function $\Phi$ is a rational function of $\bs t$, depending on
parameters $\bs z$. The master function is $\Sigma_{\bs l}$-invariant
with respect to permutations of variables with the same upper index.

The critical set of the master function with respect to variables $\bs t$
is $\Sigma_{\bs l}$-invariant.

If $\bs z$ consists of real numbers, then the critical set is
invariant with respect to complex conjugation.

\begin{conj}\label{our conj}
If $\bs z$ consists of real numbers, then every orbit of critical points is invariant
with respect to complex conjugation.
\end{conj}

For a critical point $\bs t$, define the tuple $\bs y^{\bs t} = (y_1,\dots , y_{r})$
of polynomials in variable $x$ by formulas of Section \ref{construction of spaces section}.
Conjecture \ref{our conj} can be restated as follows. If $\bs z$ consists of
real numbers
and $\bs t$ is a critical point, then the tuple $\bs y^{\bs t}$ consists
of real polynomials.

\medskip
Theorems \ref{theorem 1} and \ref{thm construction of spaces} imply this conjecture for
$\g = \slg_{r+1}$.
In the same way Theorems \ref{theorem 1} and \ref{thm construction of spaces}
imply Conjecture \ref{our conj}
for $\g$ of type $B_r$ and $C_r$, see Section 7 in \cite{MV2}.

\end{document}